\documentclass[reqno, 10pt,aps,prb, preprint]{amsart}
\usepackage{amsmath}
\usepackage{amssymb}
\newtheorem{theorem}{Theorem}[section]
\newtheorem{remark}{Remark}[section]
\newtheorem{lemma}[theorem]{Lemma}

\newtheorem{corollary}[theorem]{Corollary}

\begin{document}
\title[MHD and Navier-Stokes Systems]{Remarks on the regularity criteria of generalized MHD and Navier-Stokes systems}
\author{Kazuo Yamazaki}  
\date{}
\thanks{The author expresses gratitude to Professor Jiahong Wu and Professor David Ullrich for their teaching.}
\maketitle

\begin{abstract}
We study the regularity criteria of the three dimensional generalized MHD and Navier-Stokes systems. In particular, we show that the regularity criteria of the generalized MHD system may be reduced to depend only on two diagonal entries of the Jacobian matrix of the velocity vector field or one vorticity component and one entry of the Jacobian matrix of the velocity vector field.

\vspace{5mm}

\end{abstract}
\footnote{2000MSC : 35B65, 35Q35, 35Q86}
\footnote{Department of Mathematics, Oklahoma State University, 401 Mathematical Sciences, Stillwater, OK 74078, USA}

\section{Introduction and statement of results}

We study the generalized magnetohydrodynamics (MHD) and Navier-Stokes (NSE) systems in $\mathbb{R}^{3}$:

\begin{equation}
\begin{cases}
\partial_{t}u + (u\cdot\nabla) u -(b\cdot\nabla)b + \nabla p +  \nu\Lambda^{2\alpha} u = 0\\
\partial_{t}b + (u\cdot\nabla) b -(b\cdot\nabla)u +  \eta\Lambda^{2\beta} b = 0\\
\nabla\cdot u = \nabla\cdot b = 0, \hspace{5mm} u(x,0) = u_{0}(x), \hspace{5mm} b(x,0) = b_{0}(x)
\end{cases}
\end{equation}

\begin{equation}
\begin{cases}
\partial_{t}u + (u\cdot\nabla)u + \nabla p + \nu\Lambda^{2\alpha} u = 0\\
\nabla \cdot u = 0, \hspace{5mm} u(x,0) = u_{0}(x)
\end{cases}
\end{equation}

where $u: \mathbb{R}^{3}\times \mathbb{R}^{+} \mapsto \mathbb{R}^{3}$ represents velocity vector field, $b: \mathbb{R}^{3}\times \mathbb{R}^{+} \mapsto \mathbb{R}^{3}$ the magnetic vector field, $p: \mathbb{R}^{3} \times \mathbb{R}^{+} \mapsto \mathbb{R}$ the pressure scalar field and $\nu, \eta > 0$ the kinematic viscosity and diffusivity constants respectively. The operator $\Lambda = (-\Delta)^{\frac{1}{2}}$ is a fractional Laplacian with power $\alpha, \beta > 0$ as parameters. Without loss of generality, we set $\nu = \eta = 1$ throughout the rest of the paper. 

The global regularity issue of these systems remain one of the most challenging outstanding open problems in mathematical analysis. In two dimensional case, both MHD and NSE admit a unique global strong solution respectively; however, in three dimensional case, such results hold only locally in time (e.g. [21]).

Starting from the pioneering works of Serrin in [18] and [19] on NSE, much effort was devoted to provide sufficient conditions for a strong solution to exist globally in time and similarly for MHD (cf. e.g. [1], [2], [3], [10], [14], [15], [20], [24], [28] in case of NSE and [7], [9], [11], [22] in case of MHD). In particular, recently in [4] it was shown that the global regularity issue of the solution to NSE may depend only on one entry of the Jacobian matrix of the velocity vector field while in [5] the global regularity issue of the solution to MHD only on a partial derivative of u in $x_{3}$-direction. In relevance to such component reduction type results of regularity criteria, we mention that recently we have seen developments in the case of active scalars as well (cf. [23]).  

However, to the best of our knowledge, it is not known whether the regularity criteria of MHD system may be reduced to rely on only the entries of the Jacobian matrix of velocity vector field with number less than three. In fact, there is no regularity criteria for the system (1) even in terms of $\nabla u_{3}$, although it is known to exist for the NSE system (cf. e.g. [13]). In [12] the authors obtained partial results toward direction. 

Moreover, it is not clear whether the regularity criteria in terms of one entry in the Jacobian matrix for the classical NSE may be generalized to the case with a fractional Laplacian. Because the system (2) with $\alpha \geq \frac{5}{4}$ admits a unique global solution (cf. [21]), it is of interest if we may generalize such a result for the case $\alpha \in (1, \frac{5}{4})$. We answer these questions: 

\begin{theorem}
Let $\alpha, \beta \in (1, \frac{5}{4})$. Suppose the solution $(u,b)(x,t)$ solves (1) in $[0,T]$ and 

\begin{equation}
\int_{0}^{T}\lVert \partial_{2}u_{2}\rVert_{L^{p}}^{r} + \lVert \partial_{3}u_{3}\rVert_{L^{p}}^{r}dt < \infty
\end{equation}

for $3 \leq p < \infty$ and 

\begin{eqnarray*}
\frac{3}{p} + \frac{2\alpha}{r} &\leq& \min\{\frac{3}{p} + \frac{[5\alpha(1-\frac{1}{p}) + 4\alpha^{2}(1-\frac{1}{p}) - 10 + 4\alpha]}{2(5-2\alpha)},\\
&&\frac{3}{p} + \frac{\alpha}{\beta}\frac{[5\beta (1-\frac{1}{p}) + 4\alpha\beta (1-\frac{1}{p}) - 10 + 4\alpha]}{2(5-2\alpha)},\\
&&\frac{3}{p} + \frac{[5\alpha (1-\frac{1}{p}) + 4\alpha\beta(1-\frac{1}{p}) - 10 + 4\beta]}{2(5-2\beta)},\\
&&\frac{3}{p} + \frac{\alpha}{\beta}\frac{[5\beta(1-\frac{1}{p}) + 4\beta^{2}(1-\frac{1}{p}) - 10 + 4\beta]}{2(5-2\beta)}\},
\end{eqnarray*}

then there is no singularity up to time $T$. 
\end{theorem}

The modification of the proof of Theorem 1.1 for the classical MHD system is possible:

\begin{theorem}
Let $\alpha, \beta = 1$.  Suppose the solution $(u,b)(x,t)$ solves (1) in $[0,T]$ and 

\begin{equation}
\int_{0}^{T}\lVert \partial_{2}u_{2}\rVert_{L^{p}}^{r} + \lVert \partial_{3}u_{3}\rVert_{L^{p}}^{r}dt < \infty
\end{equation}

for $3 < p < \infty$ and 

\begin{equation*}
\frac{3}{p} + \frac{2}{r} \leq \frac{3}{2p} + \frac{1}{2},
\end{equation*}

then there is no singularity up to time $T$. 

\end{theorem}

We state an immediate interesting corollary of Theorem 1.2 which does not seem to follow from the work of [5] or [12]: 

\begin{corollary}
Let $\alpha, \beta = 1$. Suppose the solution $(u,b)(x,t)$ solves (1) in [0,T] and 

\begin{equation*}
\int_{0}^{T}\lVert w_{3}\rVert_{L^{p}}^{r} + \lVert \partial_{3}u_{3}\rVert_{L^{p}}^{r}dt < \infty
\end{equation*}

for $3 < p < \infty$ and 

\begin{equation*}
\frac{3}{p} + \frac{2}{r} \leq \frac{3}{2p} + \frac{1}{2} 
\end{equation*}

where $w_{3} = \partial_{1}u_{2} - \partial_{2}u_{1}$, then there is no singularity up to time T. 

\end{corollary} 

\begin{theorem}
Let $\alpha \in (1, \frac{5}{4})$. Suppose the solution $u(x,t)$ solves (2) in $[0,T]$ and 

\begin{equation}
\int_{0}^{T}\lVert\partial_{3}u_{3}\rVert_{L^{p}}^{r}dt < \infty
\end{equation}

where $2 < p < \infty$ and 

\begin{equation*}
\frac{3}{p} + \frac{2\alpha}{r} \leq \frac{3}{p} + \left(\frac{p-2}{p}\right)\frac{\alpha(5+4\alpha)}{4(5-2\alpha)},
\end{equation*}

then there is no singularity up to time T. 

\end{theorem}

\begin{remark}

\begin{enumerate}

\item The key to the proof of theorems of this type is an appropriate decomposition of nonlinear terms. It is not clear whether a direct extension of the proof in [4] is possible due to the complex structure of the four nonlinear terms of (1), as discussed in [12]. Our approach is based on an observation that upon $\lVert\nabla_{h}u\rVert_{L^{2}}^{2} + \lVert\nabla_{h}b\rVert_{L^{2}}^{2}$ estimate, every nonlinear term has $u$ involved. Hence, making use of the incompressibility of both $u$ and $b$, we may separate $u_{1}, u_{2}$ and $u_{3}$. Our second observation is that Lemma 2.2 below due to [4], of which originally $i=3$, may be used for $i$ any direction. Thus, we can use this lemma to concentrate the regularity dependence on $\partial_{2}u_{2}$ and $\partial_{3}u_{3}$. 
  
\item From the proof, it becomes clear that in fact we could have selected any one of the three partial derivatives of $u_{1}, u_{2}$ and $u_{3}$. Thus, for Theorem 1.2, we also proved the criteria in terms of $\partial_{3}u$ which is the result from [5]; hence, our results are more general. Moreover, Theorem 1.1 may be seen as a component reduction type result of the work of [8], [25], and [27]. Moreover, our proof may be extended to a regularity criteria of a component and a partial derivative, e.g. $u_{3}$ and $\partial_{3}u_{3}$ in the case of Theorem 1.4 as done in [26]; we chose to state the case of only partial derivatives for simplicity. 

\item The lower bound of p in the Theorem 1.1 may be optimized furthermore in terms of $\alpha$ and $\beta$; we chose to state so for simplicity. 

\item Concerning Corollary 1.3, we refer readers to [6] for similar result in the case of the NSE.

\end{enumerate}

\end{remark}

In the next sections, we list a few lemmas and thereafter prove our theorems. 

\section{Preliminaries}

We denote by $\nabla_{h}$ the horizontal gradient while $\Delta_{h}$ the horizontal Laplacian, i.e. 

\begin{equation*}
\nabla_{h} = (\partial_{1}, \partial_{2}), \hspace{5mm} \Delta_{h} = \partial_{11}^{2} + \partial_{22}^{2}
\end{equation*}

Moreover, we denote for simplicity 

\begin{equation*}
X(t) = \lVert \nabla u(t)\rVert_{L^{2}}^{2} + \lVert \nabla b(t)\rVert_{L^{2}}^{2}, \hspace{5mm} Y(t) = \lVert\nabla_{h}u(t)\rVert_{L^{2}}^{2} + \lVert \nabla_{h}b(t)\rVert_{L^{2}}^{2}
\end{equation*} 

and 

\begin{equation*}
Z(t) = \lVert \partial_{2}u_{2}(t)\rVert_{L^{s}}^{\frac{2}{\gamma -1}} + \lVert \partial_{3}u_{3}(t)\rVert_{L^{s}}^{\frac{2}{\gamma -1}}
\end{equation*}

\begin{lemma}
(cf. [14]) Let $u \in H^{2}(\mathbb{R}^{3})$ be smooth and divergence free. Then 

\begin{eqnarray*}
&&\sum_{i,j=1}^{2}\int_{\mathbb{R}^{3}}u_{i}\partial_{i}u_{j}\Delta_{h}u_{j}\\ 
&=& \sum_{i,j=1}^{2}\frac{1}{2}\int_{\mathbb{R}^{3}}\partial_{i}u_{j}\partial_{i}u_{j} \partial_{3}u_{3} - \int_{\mathbb{R}^{3}}\partial_{1}u_{1}\partial_{2}u_{2}\partial_{3}u_{3} + \int_{\mathbb{R}^{3}}\partial_{1}u_{2}\partial_{2}u_{1}\partial_{3}u_{3}
\end{eqnarray*}

\end{lemma}

\begin{lemma}
(cf. [4], [26]) For $f, g,  \in C_{c}^{\infty}(\mathbb{R}^{3})$, we have 

\begin{equation*}
\lvert \int_{\mathbb{R}^{3}}fgh\rvert \leq c\lVert f\rVert_{L^{q}}^{\frac{\gamma - 1}{\gamma}}\lVert\partial_{i}f\rVert_{L^{s}}^{\frac{1}{\gamma}}\lVert g\rVert_{L^{2}}^{\frac{\gamma -2}{\gamma}}\lVert \partial_{j}g\rVert_{L^{2}}^{\frac{1}{\gamma}}\lVert\partial_{k}g\rVert_{L^{2}}^{\frac{1}{\gamma}}\lVert h\rVert_{L^{2}}
\end{equation*}

where 

\begin{equation*}
2 < \gamma, \hspace{5mm} 1 \leq q, \hspace{5mm} s \leq \infty, \hspace{5mm} \frac{\gamma -1}{q} + \frac{1}{s} = 1
\end{equation*}

and $i, j$ and $k$ are any combinations of $1, 2$ and $3$. 
\end{lemma}

The proof of the following elementary inequality is simple and we omit it:

\begin{lemma}
For $0 \leq p < \infty$ and $a, b \geq 0$,

\begin{equation*}
(a+b)^{p} \leq 2^{p}(a^{p} + b^{p})
\end{equation*}

\end{lemma}

\section{Proofs}

\subsection{Proof of Theorem 1.1}

We start by taking an inner product of the first equation in (1) with u and the second with b and integrating in time to obtain

\begin{equation}
\sup_{t\in [0,T]}\lVert u(t)\rVert_{L^{2}}^{2} + \lVert b(t)\rVert_{L^{2}}^{2} + \int_{0}^{T}\lVert \Lambda^{\alpha}u\rVert_{L^{2}}^{2} + \lVert \Lambda^{\beta}b\rVert_{L^{2}}^{2}dt \leq c(u_{0}, b_{0})
\end{equation}

\subsubsection{Estimate of $\lVert \nabla_{h}u\rVert_{L^{2}}^{2} + \lVert \nabla_{h}b\rVert_{L^{2}}^{2}$}

Local well-posedness is shown in [21]. We devote our effort to obtain $H^{1}$ estimate below. We take an inner product of the first equation in (1) with $-\Delta_{h}u$ and the second with $-\Delta_{h}b$ to obtain 

\begin{eqnarray*}
&&\frac{1}{2}\partial_{t}Y + \lVert\Lambda^{\alpha}\nabla_{h} u\rVert_{L^{2}}^{2} + \lVert\Lambda^{\beta}\nabla_{h} b\rVert_{L^{2}}^{2}\\
&=& \int (u\cdot\nabla)u\cdot\Delta_{h}u - (b\cdot\nabla)b\cdot\Delta_{h}u + (u\cdot\nabla)b\cdot\Delta_{h}b - (b\cdot\nabla)u\cdot\Delta_{h}b = \sum_{i=1}^{4}J_{i}
\end{eqnarray*}

For $J_{1}$, we notice that applying Lemma 2.1 and integrating by parts implies

\begin{eqnarray*}
J_{1} \leq c\int\lvert u_{3}\rvert \lvert \nabla u\rvert\lvert\nabla\nabla_{h}u\rvert 
\end{eqnarray*}

For $J_{2},J_{3},J_{4}$, we decompose them as follows:

\begin{eqnarray*}
J_{2} + J_{4} &=& -\int (b\cdot\nabla)b\cdot\Delta_{h}u + (b\cdot\nabla)u\cdot\Delta_{h}b\\
&=& \sum_{i,j=1}^{3}\sum_{k=1}^{2}\int \partial_{k}b_{i}\partial_{i}b_{j}\partial_{k}u_{j} + \partial_{k}b_{i}\partial_{i}u_{j}\partial_{k}b_{j}
\end{eqnarray*}

due to the incompressibility of b. We integrate by parts once more to obtain

\begin{eqnarray*}
J_{2} + J_{4} &=&  -\sum_{i,j=1}^{3}\sum_{k=1}^{2}\int\partial_{kk}^{2}b_{i}\partial_{i}b_{j}u_{j} + \partial_{k}b_{i}\partial_{ik}^{2}b_{j}u_{j} + \partial_{ik}^{2}b_{i}u_{j} \partial_{k}b_{j} + \partial_{k}b_{i}u_{j}\partial_{ik}^{2}b_{j}\\
&\leq& c\int \lvert \nabla\nabla_{h}b\rvert \lvert \nabla b\rvert (\lvert u_{1}\rvert + \lvert u_{2}\rvert + \lvert u_{3}\rvert)
\end{eqnarray*}

Similarly, integrating by parts and using incompressibility of u, we obtain

\begin{eqnarray*}
J_{3} \leq c\int \lvert \nabla\nabla_{h}b\rvert \lvert \nabla b\rvert (\lvert u_{1}\rvert + \lvert u_{2}\rvert + \lvert u_{3}\rvert)
\end{eqnarray*}

Now we apply Lemma 2.2 with 

\begin{equation*}
f = \lvert u_{3}\rvert, g = \lvert\nabla u\rvert, h = \lvert\nabla\nabla_{h}u\rvert, i = 3, j = 1, k = 2, q = 2
\end{equation*}

to bound using (6)

\begin{eqnarray*}
J_{1} &\leq& c\lVert \partial_{3}u_{3}\rVert_{L^{s}}^{\frac{1}{\gamma}}\lVert \nabla u\rVert_{L^{2}}^{\frac{\gamma - 2}{\gamma}}\lVert \Delta u\rVert_{L^{2}}^{\frac{1}{\gamma}}\lVert \nabla\nabla_{h}u\rVert_{L^{2}}^{\frac{\gamma + 1}{\gamma}}\\
&\leq& c\lVert \partial_{3}u_{3}\rVert_{L^{s}}^{\frac{1}{\gamma}}\lVert \nabla u\rVert_{L^{2}}^{\frac{\gamma\alpha - \alpha - 1}{\gamma\alpha}}\lVert\Lambda^{\alpha}\nabla u\rVert_{L^{2}}^{\frac{1}{\gamma}\frac{1}{\alpha}}\lVert\nabla_{h}u\rVert_{L^{2}}^{(\frac{\gamma + 1}{\gamma})(1-\frac{1}{\alpha})}\lVert\Lambda^{\alpha}\nabla_{h}u\rVert_{L^{2}}^{(\frac{\gamma + 1}{\gamma})\frac{1}{\alpha}}\\
&\leq&\epsilon\lVert\Lambda^{\alpha}\nabla_{h}u\rVert_{L^{2}}^{2} + \frac{1}{4} \lVert\nabla_{h}u\rVert_{L^{2}}^{2} +c\lVert\partial_{3}u_{3}\rVert_{L^{s}}^{\frac{2}{\gamma - 1}}\lVert \nabla u\rVert_{L^{2}}^{\frac{2(\gamma\alpha - \alpha - 1)}{(\gamma - 1)\alpha}}\lVert\Lambda^{\alpha}\nabla u\rVert_{L^{2}}^{\frac{2}{(\gamma - 1)\alpha}}
\end{eqnarray*}

where we used Gagliardo-Nirenberg inequalities and Young's inequalities. Similarly, applying Lemma 2.2 and using (6) we obtain

\begin{eqnarray*}
J_{2} + J_{3} + J_{4} &\leq& c(\lVert \partial_{1}u_{1}\rVert_{L^{s}}^{\frac{1}{\gamma}}\lVert \nabla b\rVert_{L^{2}}^{\frac{\gamma -2}{\gamma}}\lVert \nabla \nabla_{h}b\rVert_{L^{2}}^{\frac{\gamma + 1}{\gamma}}\lVert\Delta b\rVert_{L^{2}}^{\frac{1}{\gamma}}\\
&&+\lVert \partial_{2}u_{2}\rVert_{L^{s}}^{\frac{1}{\gamma}}\lVert \nabla b\rVert_{L^{2}}^{\frac{\gamma -2}{\gamma}}\lVert \nabla \nabla_{h}b\rVert_{L^{2}}^{\frac{\gamma + 1}{\gamma}}\lVert\Delta b\rVert_{L^{2}}^{\frac{1}{\gamma}}\\
&&+\lVert \partial_{3}u_{3}\rVert_{L^{s}}^{\frac{1}{\gamma}}\lVert \nabla b\rVert_{L^{2}}^{\frac{\gamma -2}{\gamma}}\lVert \nabla \nabla_{h}b\rVert_{L^{2}}^{\frac{\gamma + 1}{\gamma}}\lVert\Delta b\rVert_{L^{2}}^{\frac{1}{\gamma}})
\end{eqnarray*}

Now we use divergence-free condition of u and apply Gagliardo-Nirenberg inequalities  and Young's inequalities to bound by 

\begin{eqnarray*}
&&c(\lVert\partial_{2}u_{2}\rVert_{L^{s}}^{\frac{1}{\gamma}}\lVert\nabla b\rVert_{L^{2}}^{\frac{\gamma\beta - \beta - 1}{\gamma\beta}}\lVert\nabla_{h}b\rVert_{L^{2}}^{(1+\frac{1}{\gamma})(1-\frac{1}{\beta})}\lVert\Lambda^{\beta}\nabla_{h}b\rVert_{L^{2}}^{(\frac{\gamma + 1}{\gamma})(\frac{1}{\beta})}\lVert\Lambda^{\beta}\nabla b\rVert_{L^{2}}^{\frac{1}{\gamma}\frac{1}{\beta}}\\
&&+\lVert\partial_{3}u_{3}\rVert_{L^{s}}^{\frac{1}{\gamma}}\lVert\nabla b\rVert_{L^{2}}^{\frac{\gamma\beta - \beta - 1}{\gamma\beta}}\lVert\nabla_{h}b\rVert_{L^{2}}^{(1+\frac{1}{\gamma})(1-\frac{1}{\beta})}\lVert\Lambda^{\beta}\nabla_{h}b\rVert_{L^{2}}^{(\frac{\gamma + 1}{\gamma})(\frac{1}{\beta})}\lVert\Lambda^{\beta}\nabla b\rVert_{L^{2}}^{\frac{1}{\gamma}\frac{1}{\beta}})\\
&\leq& \epsilon\lVert\Lambda^{\beta}\nabla_{h}b\rVert_{L^{2}}^{2} + \frac{1}{4}\lVert\nabla_{h}b\rVert_{L^{2}}^{2}\\
&&+c(\lVert\partial_{2}u_{2}\rVert_{L^{2}}^{\frac{2}{\gamma - 1}}\lVert\nabla b\rVert_{L^{2}}^{\frac{2(\gamma\beta - \beta - 1)}{\beta(\gamma - 1)}}\lVert\Lambda^{\beta}\nabla b\rVert_{L^{2}}^{\frac{2}{\beta(\gamma - 1)}} + \lVert\partial_{3}u_{3}\rVert_{L^{2}}^{\frac{2}{\gamma - 1}}\lVert\nabla b\rVert_{L^{2}}^{\frac{2(\gamma\beta - \beta - 1)}{\beta(\gamma - 1)}}\lVert\Lambda^{\beta}\nabla b\rVert_{L^{2}}^{\frac{2}{\beta(\gamma - 1)}})
\end{eqnarray*}

In sum, for $\epsilon > 0$ sufficiently small, 

\begin{eqnarray*}
&&\partial_{t}Y - Y + \lVert \Lambda^{\alpha}\nabla_{h}u\rVert_{L^{2}}^{2} + \lVert \Lambda^{\beta}\nabla_{h}b\rVert_{L^{2}}^{2}\\
&\leq& cZ(\lVert \nabla u\rVert_{L^{2}}^{\frac{2(\gamma\alpha - \alpha - 1)}{\alpha(\gamma - 1)}}\lVert\Lambda^{\alpha}\nabla u\rVert_{L^{2}}^{\frac{2}{\alpha(\gamma - 1)}} + \lVert\nabla b\rVert_{L^{2}}^{\frac{2(\gamma\beta - \beta - 1)}{\beta(\gamma-1)}}\lVert\Lambda^{\beta}\nabla b\rVert_{L^{2}}^{\frac{2}{\beta(\gamma - 1)}})
\end{eqnarray*}

That is,

\begin{eqnarray*}
&&\sup_{t\in [0,T]}Y(t) + \int_{0}^{T}\lVert \Lambda^{\alpha}\nabla_{h}u\rVert_{L^{2}}^{2} + \lVert\Lambda^{\beta}\nabla_{h}b\rVert_{L^{2}}^{2}dt\\
&\leq& c(T) + c\int_{0}^{T}Z(\lVert \nabla u\rVert_{L^{2}}^{\frac{2(\gamma\alpha - \alpha - 1)}{\alpha(\gamma - 1)}} \lVert \Lambda^{\alpha}\nabla u\rVert_{L^{2}}^{\frac{2}{\alpha(\gamma -1)}} + \lVert \nabla b\rVert_{L^{2}}^{\frac{2(\gamma\beta - \beta - 1)}{\beta(\gamma -1)}}\lVert\Lambda^{\beta}\nabla b\rVert_{L^{2}}^{\frac{2}{\beta(\gamma-1)}})dt
\end{eqnarray*}
 
\subsubsection{Estimate of $\lVert \nabla u\rVert_{L^{2}}^{2} + \lVert \nabla b\rVert_{L^{2}}^{2}$}

Next, we take inner products of the first equation in (1) with $-\Delta u$ and the second with $-\Delta b$ to obtain

\begin{eqnarray*}
&&\frac{1}{2}\partial_{t}X + \lVert \Lambda^{\alpha}\nabla u\rVert_{L^{2}}^{2} + \lVert \Lambda^{\beta}\nabla b\rVert_{L^{2}}^{2}\\
&=& \int (u\cdot\nabla )u\cdot\Delta_{h}u + (u\cdot\nabla)u\cdot\partial_{33}^{2}u -(b\cdot\nabla)b\cdot\Delta_{h}u - (b\cdot\nabla)b\cdot\partial_{33}^{2}u\\
&&+(u\cdot\nabla)b\cdot\Delta_{h}b + (u\cdot\nabla)b\cdot\partial_{33}^{2}b -(b\cdot\nabla)u\cdot\Delta_{h}b-(b\cdot\nabla)u\cdot\partial_{33}^{2}b
\end{eqnarray*}

We see that 

\begin{eqnarray*}
\int (u\cdot\nabla u)\cdot\partial_{33}^{2}u &=& \int u\cdot\nabla_{h}u\cdot\partial_{33}^{2}u - \frac{1}{2}(\partial_{3}u_{3})(\partial_{3}u)^{2}\\
&=& -\int\partial_{3}u\cdot\nabla_{h}u\cdot\partial_{3}u - \frac{1}{2}(\nabla_{h}\cdot u)(\partial_{3}u)^{2} - \frac{1}{2}(\partial_{1}u_{1} + \partial_{2}u_{2})(\partial_{3}u)^{2}\\
&\leq& c\int\lvert \partial_{3}u\rvert^{2}\lvert\nabla_{h}u\rvert
\end{eqnarray*}

Similarly

\begin{eqnarray*}
\int (u\cdot\nabla)b\cdot\partial_{33}^{2}b &=& \int (u\cdot\nabla_{h})b\cdot\partial_{33}^{2}b - \frac{1}{2}\partial_{3}u_{3}(\partial_{3}b)^{2}\\
&=& -\int\partial_{3}u\cdot\nabla_{h}b\cdot\partial_{3}b - \frac{1}{2}(\nabla_{h}\cdot u)(\partial_{3}b)^{2} - \frac{1}{2}(\partial_{1}u_{1} + \partial_{2}u_{2})(\partial_{3}b)^{2}\\
&\leq& c\int \lvert \partial_{3}u\rvert\lvert\nabla_{h}b\rvert \lvert \partial_{3}b\rvert + \lvert \nabla_{h}u\rvert \lvert \partial_{3}b\rvert^{2}
\end{eqnarray*}

Next, we combine two other terms:

\begin{eqnarray*}
&&\int (b\cdot\nabla)b\cdot\partial_{33}^{2}u + (b\cdot\nabla)u\cdot\partial_{33}^{2}b\\
&=& \int (b\cdot\nabla_{h})b\cdot\partial_{33}^{2}u + (b\cdot\nabla_{h})u\cdot\partial_{33}^{2}b + b_{3}\partial_{3}(\partial_{3}b\cdot\partial_{3}u)\\
&=& -\int (\partial_{3}b\cdot\nabla_{h})b\cdot\partial_{3}u + (\partial_{3}b\cdot\nabla_{h})u\cdot\partial_{3}b\\
&&-(\nabla_{h}\cdot b)(\partial_{3}b\cdot\partial_{3}u) - (\partial_{1}b_{1} + \partial_{2}b_{2})(\partial_{3}b\cdot\partial_{3}u)\\
&\leq& c\int\lvert \partial_{3}b\rvert\lvert\nabla_{h}b\rvert\lvert\partial_{3}u\rvert + \lvert \partial_{3}b\rvert^{2}\lvert \nabla_{h}u\rvert
\end{eqnarray*}

Thus, we have shown

\begin{eqnarray*}
&&\int (u\cdot\nabla)u\cdot\partial_{33}^{2}u - (b\cdot\nabla)b\cdot\partial_{33}^{2}u +(u\cdot\nabla)b\cdot\partial_{33}^{2}b - (b\cdot\nabla)u\cdot\partial_{33}^{2}b\\
&\leq& c(\lVert \nabla_{h}u\rVert_{L^{2}}\lVert\nabla u\rVert_{L^{4}}^{2} + \lVert \nabla u\rVert_{L^{4}}^{2}\lVert\nabla_{h}b\rVert_{L^{2}} + \lVert \nabla_{h}b\rVert_{L^{2}}\lVert \nabla b\rVert_{L^{4}}^{2} + \lVert \nabla_{h}u\rVert_{L^{2}}\lVert \nabla b\rVert_{L^{4}}^{2})\\
&\leq& c(\lVert \nabla_{h}u\rVert_{L^{2}}\lVert \nabla u\rVert_{L^{\frac{6}{5-2\alpha}}}^{(\frac{4\alpha - 3}{4})2}\lVert\Lambda^{\alpha}u\rVert_{L^{6}}^{(\frac{7-4\alpha}{4})2} + \lVert \nabla u\rVert_{L^{\frac{6}{5-2\alpha}}}^{(\frac{4\alpha - 3}{4})2}\lVert \Lambda^{\alpha}u\rVert_{L^{6}}^{(\frac{7-4\alpha}{4})2}\lVert\nabla_{h}b\rVert_{L^{2}}\\
&&+\lVert \nabla_{h}b\rVert_{L^{2}}\lVert \nabla b\rVert_{L^{\frac{6}{5-2\beta}}}^{(\frac{4\beta - 3}{4})2}\lVert \Lambda^{\beta}b\rVert_{L^{6}}^{(\frac{7-4\beta}{4})2} + \lVert \nabla_{h}u\rVert_{L^{2}}\lVert\nabla b\rVert_{L^{\frac{6}{5-2\beta}}}^{(\frac{4\beta - 3}{4})2}\lVert\Lambda^{\beta}b\rVert_{L^{6}}^{(\frac{7-4\beta}{4})2})\\
&\leq& c(\lVert \nabla_{h}u\rVert_{L^{2}}\lVert \Lambda^{\alpha} u\rVert_{L^{2}}^{(\frac{4\alpha - 3}{4})2}\lVert\Lambda^{\alpha}u\rVert_{L^{6}}^{(\frac{7-4\alpha}{4})2} + \lVert \Lambda^{\alpha} u\rVert_{L^{2}}^{(\frac{4\alpha - 3}{4})2}\lVert \Lambda^{\alpha}u\rVert_{L^{6}}^{(\frac{7-4\alpha}{4})2}\lVert\nabla_{h}b\rVert_{L^{2}}\\
&&+\lVert \nabla_{h}b\rVert_{L^{2}}\lVert \Lambda^{\beta} b\rVert_{L^{2}}^{(\frac{4\beta - 3}{4})2}\lVert \Lambda^{\beta}b\rVert_{L^{6}}^{(\frac{7-4\beta}{4})2} + \lVert \nabla_{h}u\rVert_{L^{2}}\lVert\Lambda^{\beta} b\rVert_{L^{2}}^{(\frac{4\beta - 3}{4})2}\lVert\Lambda^{\beta}b\rVert_{L^{6}}^{(\frac{7-4\beta}{4})2})
\end{eqnarray*}

by H$\ddot{o}$lder's inequalities, Gagliardo-Nirenberg inequalities and the Sobolev embeddings. Next, we use the well-known inequality of 

\begin{equation}
\lVert f\rVert_{L^{6}} \leq c\lVert\nabla_{h}f\rVert_{L^{2}}^{\frac{2}{3}}\lVert\partial_{3}f\rVert_{L^{2}}^{\frac{1}{3}}
\end{equation}

(cf. [5]) to bound by a constant multiples of 

\begin{eqnarray*}
&&\lVert \nabla_{h}u\rVert_{L^{2}}\lVert\Lambda^{\alpha}u\rVert_{L^{2}}^{\frac{4\alpha - 3}{2}}\lVert\Lambda^{\alpha}\nabla_{h}u\rVert_{L^{2}}^{\frac{7-4\alpha}{3}}\lVert\Lambda^{\alpha}\nabla u\rVert_{L^{2}}^{\frac{7-4\alpha}{6}}\\
&&+\lVert\nabla_{h}b\rVert_{L^{2}}\lVert\Lambda^{\alpha}u\rVert_{L^{2}}^{\frac{4\alpha - 3}{2}}\lVert\Lambda^{\alpha}\nabla_{h}u\rVert_{L^{2}}^{\frac{7-4\alpha}{3}}\lVert\Lambda^{\alpha}\nabla u\rVert_{L^{2}}^{\frac{7-4\alpha}{6}}\\
&&+\lVert\nabla_{h}b\rVert_{L^{2}}\lVert\Lambda^{\beta}b\rVert_{L^{2}}^{\frac{4\beta - 3}{2}}\lVert\Lambda^{\beta}\nabla_{h}b\rVert_{L^{2}}^{\frac{7-4\beta}{3}}\lVert\Lambda^{\beta}\nabla b\rVert_{L^{2}}^{\frac{7-4\beta}{6}}\\
&&+\lVert\nabla_{h}u\rVert_{L^{2}}\lVert\Lambda^{\beta}b\rVert_{L^{2}}^{\frac{4\beta - 3}{2}}\lVert\Lambda^{\beta}\nabla_{h}b\rVert_{L^{2}}^{\frac{7-4\beta}{3}}\lVert\Lambda^{\beta}\nabla b\rVert_{L^{2}}^{\frac{7-4\beta}{6}}
\end{eqnarray*}

We combine our previous estimate and integrate in time $[0,T]$ to obtain

\begin{eqnarray*}
&&X(T) - X(0) + \int_{0}^{T}\lVert \Lambda^{\alpha}\nabla u\rVert_{L^{2}}^{2} + \lVert \Lambda^{\beta}\nabla b\rVert_{L^{2}}^{2}dt\\
&\leq& c\int_{0}^{T}Z(\lVert \nabla u\rVert_{L^{2}}^{\frac{2(\gamma\alpha - \alpha - 1)}{\alpha(\gamma  -1)}}\lVert\Lambda^{\alpha}\nabla u\rVert_{L^{2}}^{\frac{2}{\alpha(\gamma-1)}} + \lVert \nabla b\rVert_{L^{2}}^{\frac{2(\gamma\beta - \beta - 1)}{\beta(\gamma -1)}}\lVert\Lambda^{\beta}\nabla b\rVert_{L^{2}}^{\frac{2}{\beta(\gamma - 1)}})dt\\
&&+ \int_{0}^{T}Xdt\\
&&+c \int_{0}^{T}(\lVert \nabla_{h}u\rVert_{L^{2}}\lVert\Lambda^{\alpha}u\rVert_{L^{2}}^{\frac{4\alpha - 3}{2}}\lVert\Lambda^{\alpha}\nabla_{h}u\rVert_{L^{2}}^{\frac{7-4\alpha}{3}}\lVert\Lambda^{\alpha}\nabla u\rVert_{L^{2}}^{\frac{7-4\alpha}{6}}\\
&&+\lVert\nabla_{h}b\rVert_{L^{2}}\lVert\Lambda^{\alpha}u\rVert_{L^{2}}^{\frac{4\alpha - 3}{2}}\lVert\Lambda^{\alpha}\nabla_{h}u\rVert_{L^{2}}^{\frac{7-4\alpha}{3}}\lVert\Lambda^{\alpha}\nabla u\rVert_{L^{2}}^{\frac{7-4\alpha}{6}}\\
&&+\lVert\nabla_{h}b\rVert_{L^{2}}\lVert\Lambda^{\beta}b\rVert_{L^{2}}^{\frac{4\beta - 3}{2}}\lVert\Lambda^{\beta}\nabla_{h}b\rVert_{L^{2}}^{\frac{7-4\beta}{3}}\lVert\Lambda^{\beta}\nabla b\rVert_{L^{2}}^{\frac{7-4\beta}{6}}\\
&&+\lVert\nabla_{h}u\rVert_{L^{2}}\lVert\Lambda^{\beta}b\rVert_{L^{2}}^{\frac{4\beta - 3}{2}}\lVert\Lambda^{\beta}\nabla_{h}b\rVert_{L^{2}}^{\frac{7-4\beta}{3}}\lVert\Lambda^{\beta}\nabla b\rVert_{L^{2}}^{\frac{7-4\beta}{6}})dt
\end{eqnarray*}

Now we use H$\ddot{o}$lder's inequalities and (6) to bound by a constant multiples of

\begin{eqnarray*}
&&\left(\int_{0}^{T}Z^{\frac{\alpha(\gamma -1)}{\alpha(\gamma -1) -1}}\lVert \nabla u\rVert_{L^{2}}^{2}dt\right)^{\frac{\alpha(\gamma -1) - 1}{\alpha(\gamma - 1)}}\left(\int_{0}^{T}\lVert\Lambda^{\alpha}\nabla u\rVert_{L^{2}}^{2}dt\right)^{\frac{1}{\alpha(\gamma - 1)}}\\
&& +\left(\int_{0}^{T}Z^{\frac{\beta(\gamma -1)}{\beta(\gamma -1) -1}}\lVert \nabla b\rVert_{L^{2}}^{2}dt\right)^{\frac{\beta(\gamma -1) - 1}{\beta(\gamma - 1)}}\left(\int_{0}^{T}\lVert\Lambda^{\beta}\nabla b\rVert_{L^{2}}^{2}dt\right)^{\frac{1}{\beta(\gamma - 1)}} + \int_{0}^{T}Xdt\\
&&+\sup_{t\in [0,T]}\lVert\nabla_{h}u(t)\rVert_{L^{2}}\left(\int_{0}^{T}\lVert\Lambda^{\alpha}\nabla_{h}u \rVert_{L^{2}}^{2}dt\right)^{\frac{7-4\alpha}{6}}\left(\int_{0}^{T}\lVert\Lambda^{\alpha}\nabla u\rVert_{L^{2}}^{2}dt\right)^{\frac{7-4\alpha}{12}}\\
&&+\sup_{t\in [0,T]}\lVert\nabla_{h}b(t)\rVert_{L^{2}}\left(\int_{0}^{T}\lVert\Lambda^{\alpha}\nabla_{h}u \rVert_{L^{2}}^{2}dt\right)^{\frac{7-4\alpha}{6}}\left(\int_{0}^{T}\lVert\Lambda^{\alpha}\nabla u\rVert_{L^{2}}^{2}dt\right)^{\frac{7-4\alpha}{12}}\\
&&+\sup_{t\in [0,T]}\lVert\nabla_{h}b(t)\rVert_{L^{2}}\left(\int_{0}^{T}\lVert\Lambda^{\beta} \nabla_{h}b\rVert_{L^{2}}^{2}dt\right)^{\frac{7-4\beta}{6}}\left(\int_{0}^{T}\lVert\Lambda^{\beta}\nabla b\rVert_{L^{2}}^{2}dt\right)^{\frac{7-4\beta}{12}}\\
&&+\sup_{t\in [0,T]}\lVert\nabla_{h}u(t)\rVert_{L^{2}}\left(\int_{0}^{T}\lVert\Lambda^{\beta}\nabla_{h}b \rVert_{L^{2}}^{2}dt\right)^{\frac{7-4\beta}{6}}\left(\int_{0}^{T}\lVert\Lambda^{\beta}\nabla b\rVert_{L^{2}}^{2}dt\right)^{\frac{7-4\beta}{12}}
\end{eqnarray*}

The previous estimate gives us the bound of a constant multiples of 

\begin{eqnarray*}
&&\left(\int_{0}^{T}Z^{\frac{\alpha(\gamma -1)}{\alpha(\gamma -1) -1}}\lVert \nabla u\rVert_{L^{2}}^{2}dt\right)^{\frac{\alpha(\gamma -1) - 1}{\alpha(\gamma - 1)}}\left(\int_{0}^{T}\lVert\Lambda^{\alpha}\nabla u\rVert_{L^{2}}^{2}dt\right)^{\frac{1}{\alpha(\gamma - 1)}}\\
&& +\left(\int_{0}^{T}Z^{\frac{\beta(\gamma -1)}{\beta(\gamma -1) -1}}\lVert \nabla b\rVert_{L^{2}}^{2}dt\right)^{\frac{\beta(\gamma -1) - 1}{\beta(\gamma - 1)}}\left(\int_{0}^{T}\lVert\Lambda^{\beta}\nabla b\rVert_{L^{2}}^{2}dt\right)^{\frac{1}{\beta(\gamma - 1)}} + \int_{0}^{T}Xdt\\
&&+\left(\int_{0}^{T}Z(\lVert \nabla u\rVert_{L^{2}}^{\frac{2(\gamma\alpha - \alpha - 1)}{\alpha(\gamma - 1)}} \lVert \Lambda^{\alpha}\nabla u\rVert_{L^{2}}^{\frac{2}{\alpha(\gamma -1)}} + \lVert \nabla b\rVert_{L^{2}}^{\frac{2(\gamma\beta - \beta - 1)}{\beta(\gamma -1)}}\lVert\Lambda^{\beta}\nabla b\rVert_{L^{2}}^{\frac{2}{\beta(\gamma-1)}})dt\right)^{\frac{5-2\alpha}{3}}\\
&&\times\left(\int_{0}^{T}\lVert\Lambda^{\alpha}\nabla u\rVert_{L^{2}}^{2}dt\right)^{\frac{7-4\alpha}{12}}\\
&&+\left(\int_{0}^{T}Z(\lVert \nabla u\rVert_{L^{2}}^{\frac{2(\gamma\alpha - \alpha - 1)}{\alpha(\gamma - 1)}} \lVert \Lambda^{\alpha}\nabla u\rVert_{L^{2}}^{\frac{2}{\alpha(\gamma -1)}} + \lVert \nabla b\rVert_{L^{2}}^{\frac{2(\gamma\beta - \beta - 1)}{\beta(\gamma -1)}}\lVert\Lambda^{\beta}\nabla b\rVert_{L^{2}}^{\frac{2}{\beta(\gamma-1)}})dt\right)^{\frac{5-2\beta}{3}}\\
&&\times\left(\int_{0}^{T}\lVert\Lambda^{\beta}\nabla b\rVert_{L^{2}}^{2}dt\right)^{\frac{7-4\beta}{12}}
\end{eqnarray*}

On the first two terms, we use Young's inequalities to bound by 

\begin{eqnarray*}
&&c\int_{0}^{T}Z^{\frac{\alpha(\gamma - 1)}{\alpha(\gamma - 1) - 1}}\lVert \nabla u\rVert_{L^{2}}^{2}dt + c\int_{0}^{T}Z^{\frac{\beta(\gamma - 1)}{\beta(\gamma - 1) - 1}}\lVert \nabla b\rVert_{L^{2}}^{2}dt\\
&&+\frac{1}{2}\left(\int_{0}^{T}\lVert\Lambda^{\alpha}\nabla u\rVert_{L^{2}}^{2}dt + \int_{0}^{T}\lVert\Lambda^{\beta}\nabla b\rVert_{L^{2}}^{2}dt\right)
\end{eqnarray*}

On the last two terms, by H$\ddot{o}$lder's inequalities we bound by 

\begin{eqnarray*}
&&c[\left(\int_{0}^{T}Z^{\frac{\alpha(\gamma-1)}{\alpha(\gamma - 1) - 1}}\lVert\nabla u\rVert_{L^{2}}^{2}dt\right)^{\frac{\alpha(\gamma - 1) - 1}{\alpha(\gamma - 1)}}\left(\int_{0}^{T}\lVert\Lambda^{\alpha}\nabla u\rVert_{L^{2}}^{2}dt\right)^{\frac{1}{\alpha(\gamma - 1)}}\\
&&+ \left(\int_{0}^{T}Z^{\frac{\beta(\gamma-1)}{\beta(\gamma - 1) - 1}}\lVert\nabla b\rVert_{L^{2}}^{2}dt\right)^{\frac{\beta(\gamma - 1) - 1}{\beta(\gamma-1)}}\left(\int_{0}^{T}\lVert\Lambda^{\beta}\nabla b\rVert_{L^{2}}^{2}dt\right)^{\frac{1}{\beta(\gamma-1)}}]^{\frac{5-2\alpha}{3}}\\
&&\times \left(\int_{0}^{T}\lVert\Lambda^{\alpha}\nabla u\rVert_{L^{2}}^{2}dt\right)^{\frac{7-4\alpha}{12}}\\
&&+c[\left(\int_{0}^{T}Z^{\frac{\alpha(\gamma-1)}{\alpha(\gamma - 1) - 1}}\lVert\nabla u\rVert_{L^{2}}^{2}dt\right)^{\frac{\alpha(\gamma - 1) - 1}{\alpha(\gamma - 1)}}\left(\int_{0}^{T}\lVert\Lambda^{\alpha}\nabla u\rVert_{L^{2}}^{2}dt\right)^{\frac{1}{\alpha(\gamma - 1)}}\\
&&+ \left(\int_{0}^{T}Z^{\frac{\beta(\gamma-1)}{\beta(\gamma - 1) - 1}}\lVert\nabla b\rVert_{L^{2}}^{2}dt\right)^{\frac{\beta(\gamma - 1) - 1}{\beta(\gamma-1)}}\left(\int_{0}^{T}\lVert\Lambda^{\beta}\nabla b\rVert_{L^{2}}^{2}dt\right)^{\frac{1}{\beta(\gamma-1)}}]^{\frac{5-2\beta}{3}}\\
&&\times \left(\int_{0}^{T}\lVert\Lambda^{\beta}\nabla b\rVert_{L^{2}}^{2}dt\right)^{\frac{7-4\beta}{12}}
\end{eqnarray*}

Next, by Lemma 2.3 we bound by a constant multiples of 

\begin{eqnarray*}
&&\left(\int_{0}^{T}Z^{\frac{\alpha(\gamma - 1)}{\alpha(\gamma - 1) - 1}}\lVert \nabla u\rVert_{L^{2}}^{2}dt\right)^{[\frac{\alpha(\gamma - 1) - 1}{\alpha(\gamma - 1)}](\frac{5-2\alpha}{3})}\left(\int_{0}^{T}\lVert\Lambda^{\alpha}\nabla u\rVert_{L^{2}}^{2}dt\right)^{[\frac{1}{\alpha(\gamma-1)}](\frac{5-2\alpha}{3}) + \frac{7-4\alpha}{12}}\\
&&+ \left(\int_{0}^{T}Z^{\frac{\beta(\gamma - 1)}{\beta(\gamma - 1) - 1}}\lVert \nabla b\rVert_{L^{2}}^{2}dt\right)^{[\frac{\beta(\gamma - 1) - 1}{\beta(\gamma - 1)}](\frac{5-2\alpha}{3})}\left(\int_{0}^{T}\lVert\Lambda^{\beta}\nabla b\rVert_{L^{2}}^{2}dt\right)^{[\frac{1}{\beta(\gamma-1)}](\frac{5-2\alpha}{3})}\\
&&\times\left(\int_{0}^{T}\lVert\Lambda^{\alpha}\nabla u\rVert_{L^{2}}^{2}dt\right)^{\frac{7-4\alpha}{12}}\\
&&+ \left(\int_{0}^{T}Z^{\frac{\alpha(\gamma - 1)}{\alpha(\gamma - 1) - 1}}\lVert \nabla u\rVert_{L^{2}}^{2}dt\right)^{[\frac{\alpha(\gamma - 1) - 1}{\alpha(\gamma - 1)}](\frac{5-2\beta}{3})}\left(\int_{0}^{T}\lVert\Lambda^{\alpha}\nabla u\rVert_{L^{2}}^{2}dt\right)^{[\frac{1}{\alpha(\gamma-1)}](\frac{5-2\beta}{3})}\\
&&\times\left(\int_{0}^{T}\lVert\Lambda^{\beta}\nabla b\rVert_{L^{2}}^{2}dt\right)^{\frac{7-4\beta}{12}}\\
&&+\left(\int_{0}^{T}Z^{\frac{\beta(\gamma - 1)}{\beta(\gamma - 1) - 1}}\lVert \nabla b\rVert_{L^{2}}^{2}dt\right)^{[\frac{\beta(\gamma - 1) - 1}{\beta(\gamma - 1)}](\frac{5-2\beta}{3})}\left(\int_{0}^{T}\lVert\Lambda^{\beta}\nabla b\rVert_{L^{2}}^{2}dt\right)^{[\frac{1}{\beta(\gamma-1)}](\frac{5-2\beta}{3}) + \frac{7-4\beta}{12}}
\end{eqnarray*}

Now we use Young's inequalities to bound by a constant multiples of 

\begin{eqnarray*}
&&\left(\int_{0}^{T}Z^{\frac{\alpha(\gamma - 1)}{\alpha(\gamma -1) - 1}}\lVert \nabla u\rVert_{L^{2}}^{2}dt\right)^{\frac{[\alpha(\gamma -1) - 1][5-2\alpha]}{3\alpha(\gamma -1) - 5 + 2\alpha}}\left(\int_{0}^{T}\lVert\Lambda^{\alpha}\nabla u\rVert_{L^{2}}^{2}dt\right)^{(\frac{7-4\alpha}{12})\frac{3\alpha(\gamma-1)}{[3\alpha(\gamma -1) - 5 + 2\alpha]}}\\
&&+\left(\int_{0}^{T}Z^{\frac{\beta(\gamma - 1)}{\beta(\gamma -1) - 1}}\lVert \nabla b\rVert_{L^{2}}^{2}dt\right)^{\frac{[\beta(\gamma -1) - 1][5-2\alpha]}{3\beta(\gamma -1) - 5 + 2\alpha}}\left(\int_{0}^{T}\lVert\Lambda^{\alpha}\nabla u\rVert_{L^{2}}^{2}dt\right)^{(\frac{7-4\alpha}{12})\frac{3\beta(\gamma-1)}{[3\beta(\gamma -1) - 5 + 2\alpha]}}\\
&&+\left(\int_{0}^{T}Z^{\frac{\alpha(\gamma - 1)}{\alpha(\gamma -1) - 1}}\lVert \nabla u\rVert_{L^{2}}^{2}dt\right)^{\frac{[\alpha(\gamma -1) - 1][5-2\beta]}{3\alpha(\gamma -1) - 5 + 2\beta}}\left(\int_{0}^{T}\lVert\Lambda^{\beta}\nabla b\rVert_{L^{2}}^{2}dt\right)^{(\frac{7-4\beta}{12})\frac{3\alpha(\gamma-1)}{[3\alpha(\gamma -1) - 5 + 2\beta]}}\\
&&+\left(\int_{0}^{T}Z^{\frac{\beta(\gamma - 1)}{\beta(\gamma -1) - 1}}\lVert \nabla b\rVert_{L^{2}}^{2}dt\right)^{\frac{[\beta(\gamma -1) - 1][5-2\beta]}{3\beta(\gamma -1) - 5 + 2\beta}}\left(\int_{0}^{T}\lVert\Lambda^{\beta}\nabla b\rVert_{L^{2}}^{2}dt\right)^{(\frac{7-4\beta}{12})\frac{3\beta(\gamma-1)}{[3\beta(\gamma -1) - 5 + 2\beta]}}\\
&&+\int_{0}^{T}\lVert\Lambda^{\alpha}\nabla u\rVert_{L^{2}}^{2} + \lVert\Lambda^{\beta}\nabla b\rVert_{L^{2}}^{2}dt
\end{eqnarray*}

Next we will use the following Young's inequalities with 

\begin{equation*}
\begin{cases}
(\frac{7-4\alpha}{4})[\frac{\alpha(\gamma - 1)}{3\alpha(\gamma - 1) - 5 + 2\alpha}] + \frac{5\alpha(\gamma -1) + 4\alpha^{2}(\gamma -1) - 20 + 8\alpha}{4[3\alpha(\gamma -1) - 5 + 2\alpha]} = 1\\
(\frac{7-4\alpha}{4})[\frac{\beta(\gamma - 1)}{3\beta(\gamma - 1) - 5 + 2\alpha}] + \frac{5\beta(\gamma -1) + 4\alpha\beta(\gamma -1) - 20 + 8\alpha}{4[3\beta(\gamma -1) - 5 + 2\alpha]} = 1\\
(\frac{7-4\beta}{4})[\frac{\alpha(\gamma - 1)}{3\alpha(\gamma - 1) - 5 + 2\beta}] + \frac{5\alpha(\gamma -1) + 4\alpha\beta(\gamma -1) - 20 + 8\beta}{4[3\alpha(\gamma -1) - 5 + 2\beta]} = 1\\
(\frac{7-4\beta}{4})[\frac{\beta(\gamma - 1)}{3\beta(\gamma - 1) - 5 + 2\beta}] + \frac{5\beta(\gamma -1) + 4\beta^{2}(\gamma -1) - 20 + 8\beta}{4[3\beta(\gamma -1) - 5 + 2\beta]} = 1
\end{cases}
\end{equation*}

Justification of these Young's inequalities require in addition to the hypothesis of Lemma 2.2,  

\begin{equation*}
\gamma > \max\{2, \frac{20-3\alpha + 4\alpha^{2}}{5\alpha + 4\alpha^{2}}, \frac{20+ 5\beta - 8\alpha +  4\alpha\beta}{5\beta + 4\alpha\beta}, \frac{20+5\alpha - 8\beta + 4\alpha\beta}{5\alpha +4\alpha\beta}, \frac{20-3\beta+4\beta^{2}}{5\beta + 4\beta^{2}}\}
\end{equation*}

Thus, letting 

\begin{equation*}
\begin{cases}
\frac{1}{p_{1}} = \frac{4[\alpha(\gamma -1) -1][5-2\alpha]}{5\alpha(\gamma-1) + 4\alpha^{2}(\gamma -1) - 20 + 8\alpha}\\
\frac{1}{p_{2}} = \frac{4[\beta(\gamma -1) -1][5-2\alpha]}{5\beta(\gamma-1) + 4\alpha\beta(\gamma -1) - 20 + 8\alpha}\\
\frac{1}{p_{3}} = \frac{4[\alpha(\gamma -1) -1][5-2\beta]}{5\alpha(\gamma-1) + 4\alpha\beta(\gamma -1) - 20 + 8\beta}\\
\frac{1}{p_{4}} = \frac{4[\beta(\gamma -1) -1][5-2\beta]}{5\beta(\gamma-1) + 4\beta^{2}(\gamma -1) - 20 + 8\beta}
\end{cases}
\end{equation*}

we have the bound by a constant multiples of

\begin{eqnarray*}
&&\left(\int_{0}^{T}Z^{\frac{\alpha(\gamma -1)}{\alpha(\gamma -1) - 1}}\lVert \nabla u\rVert_{L^{2}}^{2p_{1}}\lVert \nabla u\rVert_{L^{2}}^{2(1-p_{1})}dt\right)^{\frac{1}{p_{1}}}\\
&&+\left(\int_{0}^{T}Z^{\frac{\beta(\gamma -1)}{\beta(\gamma -1) - 1}}\lVert \nabla b\rVert_{L^{2}}^{2p_{2}}\lVert \nabla b\rVert_{L^{2}}^{2(1-p_{2})}dt\right)^{\frac{1}{p_{2}}}\\
&&+\left(\int_{0}^{T}Z^{\frac{\alpha(\gamma -1)}{\alpha(\gamma -1) - 1}}\lVert \nabla u\rVert_{L^{2}}^{2p_{3}}\lVert \nabla u\rVert_{L^{2}}^{2(1-p_{3})}dt\right)^{\frac{1}{p_{3}}}\\
&&+\left(\int_{0}^{T}Z^{\frac{\beta(\gamma -1)}{\beta(\gamma -1) - 1}}\lVert \nabla b\rVert_{L^{2}}^{2p_{4}}\lVert \nabla b\rVert_{L^{2}}^{2(1-p_{4})}dt\right)^{\frac{1}{p_{4}}}
\end{eqnarray*}

By H$\ddot{o}$lder's inequalities we further bound by a constant multiples of

\begin{eqnarray*}
&&\left(\int_{0}^{T}Z^{\frac{\alpha(\gamma-1)}{\alpha(\gamma-1) - 1}(\frac{1}{p_{1}})}\lVert \nabla u\rVert_{L^{2}}^{2}dt\right)\left(\int_{0}^{T}\lVert \nabla u\rVert_{L^{2}}^{2}dt\right)^{\frac{1-p_{1}}{p_{1}}}\\
&&+\left(\int_{0}^{T}Z^{\frac{\beta(\gamma-1)}{\beta(\gamma-1) - 1}(\frac{1}{p_{2}})}\lVert \nabla b\rVert_{L^{2}}^{2}dt\right)\left(\int_{0}^{T}\lVert \nabla b\rVert_{L^{2}}^{2}dt\right)^{\frac{1-p_{2}}{p_{2}}}\\
&&+\left(\int_{0}^{T}Z^{\frac{\alpha(\gamma-1)}{\alpha(\gamma-1) - 1}(\frac{1}{p_{3}})}\lVert \nabla u\rVert_{L^{2}}^{2}dt\right)\left(\int_{0}^{T}\lVert \nabla u\rVert_{L^{2}}^{2}\right)^{\frac{1-p_{3}}{p_{3}}}\\
&&+\left(\int_{0}^{T}Z^{\frac{\beta(\gamma-1)}{\beta(\gamma-1) - 1}(\frac{1}{p_{4}})}\lVert \nabla b\rVert_{L^{2}}^{2}dt\right)\left(\int_{0}^{T}\lVert \nabla b\rVert_{L^{2}}^{2}dt\right)^{\frac{1-p_{4}}{p_{4}}}
\end{eqnarray*}

By Gagliardo-Nirenberg inequalities and (6), we have obtained

\begin{eqnarray*}
&&X(T) + \int_{0}^{T}\lVert \Lambda^{\alpha}\nabla u\rVert_{L^{2}}^{2} + \lVert \Lambda^{\beta}\nabla b\rVert_{L^{2}}^{2}dt\\
&\leq& X(0) + c\int_{0}^{T}Z^{\frac{\alpha(\gamma-1)}{\alpha(\gamma-1) - 1}}\lVert \nabla u\rVert_{L^{2}}^{2} + Z^{\frac{\beta(\gamma-1)}{\beta(\gamma-1)-1}}\lVert \nabla b\rVert_{L^{2}}^{2} + X\\
&&+Z^{\frac{\alpha(\gamma-1)}{\alpha(\gamma-1) - 1}(\frac{1}{p_{1}})}\lVert \nabla u\rVert_{L^{2}}^{2} + Z^{\frac{\beta(\gamma-1)}{\beta(\gamma-1)-1}(\frac{1}{p_{2}})}\lVert \nabla b\rVert_{L^{2}}^{2} \\
&&+Z^{\frac{\alpha(\gamma-1)}{\alpha(\gamma-1) - 1}(\frac{1}{p_{3}})}\lVert \nabla u\rVert_{L^{2}}^{2} + Z^{\frac{\beta(\gamma-1)}{\beta(\gamma-1)-1}(\frac{1}{p_{4}})}\lVert \nabla b\rVert_{L^{2}}^{2}dt
\end{eqnarray*}

By Lemma 2.3

\begin{eqnarray*}
&&Z^{\frac{\alpha(\gamma-1)}{\alpha(\gamma-1) - 1}(\frac{1}{p_{1}})}+ Z^{\frac{\beta(\gamma-1)}{\beta(\gamma-1)-1}(\frac{1}{p_{2}})}+Z^{\frac{\alpha(\gamma-1)}{\alpha(\gamma-1) - 1}(\frac{1}{p_{3}})}+ Z^{\frac{\beta(\gamma-1)}{\beta(\gamma-1)-1}(\frac{1}{p_{4}})}\\
&\leq& c(\lVert \partial_{2}u_{2}\rVert_{L^{s}}^{\frac{8\alpha(5-2\alpha)}{5\alpha(\gamma-1) + 4\alpha^{2}(\gamma-1) - 20 + 8\alpha}} + \lVert \partial_{3}u_{3}\rVert_{L^{s}}^{\frac{8\alpha(5-2\alpha)}{5\alpha(\gamma-1) + 4\alpha^{2}(\gamma-1) - 20 + 8\alpha}}\\
&&+\lVert \partial_{2}u_{2}\rVert_{L^{s}}^{\frac{8\beta(5-2\alpha)}{5\beta(\gamma-1) + 4\alpha\beta(\gamma-1) - 20 + 8\alpha}} + \lVert \partial_{3}u_{3}\rVert_{L^{s}}^{\frac{8\beta(5-2\alpha)}{5\beta(\gamma-1) + 4\alpha\beta(\gamma-1) - 20 + 8\alpha}}\\
&&+\lVert \partial_{2}u_{2}\rVert_{L^{s}}^{\frac{8\alpha(5-2\beta)}{5\alpha(\gamma-1) + 4\alpha\beta(\gamma-1) - 20 + 8\beta}} + \lVert \partial_{3}u_{3}\rVert_{L^{s}}^{\frac{8\alpha(5-2\beta)}{5\alpha(\gamma-1) + 4\alpha\beta(\gamma-1) - 20 + 8\beta}}\\
&&+\lVert \partial_{2}u_{2}\rVert_{L^{s}}^{\frac{8\beta(5-2\beta)}{5\beta(\gamma-1) +4\beta^{2}(\gamma-1) - 20 + 8\beta}} + \lVert \partial_{3}u_{3}\rVert_{L^{s}}^{\frac{8\beta(5-2\beta)}{5\beta(\gamma-1) + 4\beta^{2}(\gamma-1) - 20 + 8\beta}})
\end{eqnarray*}

Gronwall's inequality and (3) imply the desired result. Lastly, considering the range of $\alpha, \beta$, all the conditions of 

\begin{eqnarray*}
\gamma &>& \max\{2, \frac{20-3\alpha + 4\alpha^{2}}{5\alpha + 4\alpha^{2}}, \frac{20 + 5\beta - 8\alpha + 4\alpha\beta}{5\beta + 4\alpha\beta}, \frac{20+5\alpha - 8\beta  +4\alpha\beta}{5\alpha + 4\alpha\beta},\\
&&\frac{20-3\beta +4\beta^{2}}{5\beta + 4\beta^{2}}\}
\end{eqnarray*}

may be simplified to say that $\gamma > \frac{7}{3}$ suffices. 

\subsection{Proof of Theorem 1.2}

The proof is similar to that of Theorem 1.1 and hence we only sketch it. We start by taking an inner product of the first equation in (1) with u and the second with b and integrating in time to obtain 

\begin{equation}
\sup_{t\in [0,T]}\lVert u(t)\rVert_{L^{2}}^{2} + \lVert b(t)\rVert_{L^{2}}^{2} + \int_{0}^{T}Xdt \leq c(u_{0}, b_{0})
\end{equation}

\subsubsection{Estimate of $\lVert \nabla_{h}u\rVert_{L^{2}}^{2} + \lVert \nabla_{h}b\rVert_{L^{2}}^{2}$}

Local well-posedness is shown in [17]. We take an inner product of the first equation in (1) with $-\Delta_{h}u$ and the second with $-\Delta_{h}b$ to obtain 

\begin{eqnarray*}
&&\frac{1}{2}\partial_{t}Y + \lVert\nabla_{h}\nabla u\rVert_{L^{2}}^{2} + \lVert\nabla_{h}\nabla b\rVert_{L^{2}}^{2}\\
&=& \int (u\cdot\nabla)u\cdot\Delta_{h}u - (b\cdot\nabla)b\cdot\Delta_{h}u + (u\cdot\nabla)b\cdot\Delta_{h}b - (b\cdot\nabla)u\cdot\Delta_{h}b = \sum_{i=1}^{4}J_{i}
\end{eqnarray*}

Similarly as before, Lemma 2.2, (8) and Young's inequality give

\begin{eqnarray*}
J_{1} \leq \epsilon\lVert \nabla \nabla_{h}u\rVert_{L^{2}}^{2} + c\lVert\partial_{3}u\rVert_{L^{s}}^{\frac{2}{\gamma - 2}}\lVert\nabla u\rVert_{L^{2}}^{2}
\end{eqnarray*}

and

\begin{eqnarray*}
J_{2} + J_{3} + J_{4} &\leq& \epsilon\lVert \nabla\nabla_{h}b\rVert_{L^{2}}^{2} + cZ\lVert \nabla b\rVert_{L^{2}}^{2(\frac{\gamma - 2}{\gamma -1})}\lVert \Delta b\rVert_{L^{2}}^{\frac{2}{\gamma -1}}
\end{eqnarray*}

Thus, for $\epsilon > 0$ sufficiently small, due to the divergence-free property of $u$, 

\begin{eqnarray*}
&&\sup_{t\in [0,T]}Y(t) + \int_{0}^{T}\lVert \nabla \nabla_{h}u\rVert_{L^{2}}^{2}  +\lVert \nabla \nabla_{h}b\rVert_{L^{2}}^{2}dt\\
&\leq& Y(0) + c\int_{0}^{T}\lVert \partial_{3}u_{3}\rVert_{L^{s}}^{\frac{2}{\gamma - 2}}X +cZ\lVert \nabla b\rVert_{L^{2}}^{2(\frac{\gamma - 2}{\gamma - 1})}\lVert\Delta b\rVert_{L^{2}}^{\frac{2}{\gamma - 1}}dt  
\end{eqnarray*}

\subsubsection{Estimate of $\lVert \nabla u\rVert_{L^{2}}^{2} + \lVert \nabla b\rVert_{L^{2}}^{2}$}

Next, we take inner products of the first equation in (1) with $-\Delta u$ and the second with $-\Delta b$ to obtain

\begin{eqnarray*}
&&\frac{1}{2}\partial_{t}X + \lVert \Delta u\rVert_{L^{2}}^{2} + \lVert \Delta b\rVert_{L^{2}}^{2}\\
&=& \int (u\cdot\nabla )u\cdot\Delta_{h}u + (u\cdot\nabla)u\cdot\partial_{33}^{2}u -(b\cdot\nabla)b\cdot\Delta_{h}u - (b\cdot\nabla)b\cdot\partial_{33}^{2}u\\
&&+(u\cdot\nabla)b\cdot\Delta_{h}b + (u\cdot\nabla)b\cdot\partial_{33}^{2}b -(b\cdot\nabla)u\cdot\Delta_{h}b-(b\cdot\nabla)u\cdot\partial_{33}^{2}b
\end{eqnarray*}

Similarly as before, we can obtain

\begin{eqnarray*}
&&\int (u\cdot\nabla)u\cdot\partial_{33}^{2}u - (b\cdot\nabla)b\cdot\partial_{33}^{2}u +(u\cdot\nabla)b\cdot\partial_{33}^{2}b - (b\cdot\nabla)u\cdot\partial_{33}^{2}b\\
&\leq& c\int\lvert \partial_{3}u\rvert^{2}\lvert \nabla_{h}u\rvert + \lvert \partial_{3}u\rvert \lvert \nabla_{h}b\rvert \lvert \partial_{3}b\rvert + \lvert \nabla_{h}u\rvert\lvert\partial_{3}b\rvert^{2}
\end{eqnarray*}

Then, H$\ddot{o}$lder's inequality, (7) and the previous estimate give after integrating in time

\begin{eqnarray*}
&&X(T) + \int_{0}^{T}\lVert \Delta u\rVert_{L^{2}}^{2} + \lVert \Delta b\rVert_{L^{2}}^{2}dt\\
&\leq& X(0) + c\int_{0}^{T}\lVert\partial_{3}u_{3}\rVert_{L^{s}}^{\frac{2}{\gamma - 2}}Xdt +c\int_{0}^{T}Z\lVert\nabla b\rVert_{L^{2}}^{2(\frac{\gamma - 2}{\gamma - 1})}\lVert \Delta b \rVert_{L^{2}}^{\frac{2}{\gamma - 1}}dt\\
&&+c\int_{0}^{T}\lVert \nabla_{h}u\rVert_{L^{2}}\lVert \nabla u\rVert_{L^{2}}^{\frac{1}{2}}\lVert \nabla\nabla_{h}u\rVert_{L^{2}}\lVert\Delta u\rVert_{L^{2}}^{\frac{1}{2}}  +\lVert \nabla_{h}b\rVert_{L^{2}}\lVert\nabla u\rVert_{L^{2}}^{\frac{1}{2}}\lVert \nabla\nabla_{h}u\rVert_{L^{2}}\lVert \Delta u\rVert_{L^{2}}^{\frac{1}{2}}\\
&&+\lVert\nabla_{h}b\rVert_{L^{2}}\lVert\nabla b\rVert_{L^{2}}^{\frac{1}{2}}\lVert\nabla\nabla_{h}b\rVert_{L^{2}}\lVert\Delta b\rVert_{L^{2}}^{\frac{1}{2}} +\lVert \nabla_{h}u\rVert_{L^{2}}\lVert \nabla b\rVert_{L^{2}}^{\frac{1}{2}}\lVert\nabla\nabla_{h}b\rVert_{L^{2}}\lVert\Delta b\rVert_{L^{2}}^{\frac{1}{2}}dt
\end{eqnarray*}

Another H$\ddot{o}$lder's inequalities imply

\begin{eqnarray*}
&&X(T) + \int_{0}^{T}\lVert \Delta u\rVert_{L^{2}}^{2} + \lVert \Delta b\rVert_{L^{2}}^{2}dt\\
&\leq& X(0) + c\int_{0}^{T}\lVert\partial_{3}u_{3}\rVert_{L^{s}}^{\frac{2}{\gamma - 2}}Xdt +c\left(\int_{0}^{T}Z^{\frac{\gamma - 1}{\gamma - 2}}\lVert \nabla b\rVert_{L^{2}}^{2}dt\right)^{\frac{\gamma - 2}{\gamma - 1}}\left(\int_{0}^{T}\lVert \Delta b\rVert_{L^{2}}^{2}dt\right)^{\frac{1}{\gamma - 1}}\\
&&+c[\sup_{t}\lVert \nabla_{h}u(t)\rVert_{L^{2}}\left(\int_{0}^{T}\lVert\nabla u\rVert_{L^{2}}^{2}dt\right)^{\frac{1}{4}}\left(\int_{0}^{T}\lVert\nabla\nabla_{h}u\rVert_{L^{2}}^{2}dt\right)^{\frac{1}{2}}\left(\int_{0}^{T}\lVert\Delta u\rVert_{L^{2}}^{2}dt\right)^{\frac{1}{4}}\\
&&+\sup_{t}\lVert \nabla_{h}b(t)\rVert_{L^{2}}\left(\int_{0}^{T}\lVert\nabla u\rVert_{L^{2}}^{2}dt\right)^{\frac{1}{4}}\left(\int_{0}^{T}\lVert\nabla\nabla_{h}u\rVert_{L^{2}}^{2}dt\right)^{\frac{1}{2}}\left(\int_{0}^{T}\lVert\Delta u\rVert_{L^{2}}^{2}dt\right)^{\frac{1}{4}}\\
&&+\sup_{t}\lVert \nabla_{h}b(t)\rVert_{L^{2}}\left(\int_{0}^{T}\lVert\nabla b\rVert_{L^{2}}^{2}dt\right)^{\frac{1}{4}}\left(\int_{0}^{T}\lVert\nabla\nabla_{h}b\rVert_{L^{2}}^{2}dt\right)^{\frac{1}{2}}\left(\int_{0}^{T}\lVert\Delta b\rVert_{L^{2}}^{2}dt\right)^{\frac{1}{4}}]\\
&&+\sup_{t}\lVert \nabla_{h}u(t)\rVert_{L^{2}}\left(\int_{0}^{T}\lVert\nabla b\rVert_{L^{2}}^{2}dt\right)^{\frac{1}{4}}\left(\int_{0}^{T}\lVert\nabla\nabla_{h}b\rVert_{L^{2}}^{2}dt\right)^{\frac{1}{2}}\left(\int_{0}^{T}\lVert\Delta b\rVert_{L^{2}}^{2}dt\right)^{\frac{1}{4}}]
\end{eqnarray*}

Next, (8), Young's inequality and previous estimates imply

\begin{eqnarray*}
&&X(T) + \int_{0}^{T}\lVert \Delta u\rVert_{L^{2}}^{2} + \lVert \Delta b\rVert_{L^{2}}^{2}dt\\
&\leq& X(0) + c\int_{0}^{T}\lVert\partial_{3}u_{3}\rVert_{L^{s}}^{\frac{2}{\gamma - 2}}Xdt +c\int_{0}^{T}Z^{\frac{\gamma - 1}{\gamma - 2}}\lVert \nabla b\rVert_{L^{2}}^{2}dt + \epsilon\int_{0}^{T}\lVert \Delta b\rVert_{L^{2}}^{2}dt\\
&&+c [\int_{0}^{T}\lVert \partial_{3}u_{3}\rVert_{L^{s}}^{\frac{2}{\gamma - 2}}Xdt + \int_{0}^{T}Z\lVert \nabla b\rVert_{L^{2}}^{2(\frac{\gamma - 2}{\gamma - 1})}\lVert\Delta b\rVert_{L^{2}}^{\frac{2}{\gamma - 1}}dt]\\
&&\times[\left(\int_{0}^{T}\lVert\Delta u\rVert_{L^{2}}^{2}dt\right)^{\frac{1}{4}} + \left(\int_{0}^{T}\lVert\Delta b\rVert_{L^{2}}^{2}dt\right)^{\frac{1}{4}}]\\
&=& X(0) + c\int_{0}^{T}\lVert\partial_{3}u_{3}\rVert_{L^{s}}^{\frac{2}{\gamma - 2}}Xdt +c\int_{0}^{T}Z^{\frac{\gamma - 1}{\gamma - 2}}\lVert \nabla b\rVert_{L^{2}}^{2}dt + \epsilon\int_{0}^{T}\lVert \Delta b\rVert_{L^{2}}^{2}dt\\
&&+c\int_{0}^{T}\lVert \partial_{3}u_{3}\rVert_{L^{s}}^{\frac{2}{\gamma - 2}}Xdt[\left(\int_{0}^{T}\lVert\Delta u\rVert_{L^{2}}^{2}dt\right)^{\frac{1}{4}} + \left(\int_{0}^{T}\lVert\Delta b\rVert_{L^{2}}^{2}dt\right)^{\frac{1}{4}}]\\
&&+c\int_{0}^{T}Z\lVert \nabla b\rVert_{L^{2}}^{2(\frac{\gamma - 2}{\gamma - 1})}\lVert\Delta b\rVert_{L^{2}}^{\frac{2}{\gamma - 1}}dt[\left(\int_{0}^{T}\lVert\Delta u\rVert_{L^{2}}^{2}dt\right)^{\frac{1}{4}} +\left(\int_{0}^{T}\lVert\Delta b\rVert_{L^{2}}^{2}dt\right)^{\frac{1}{4}}]
\end{eqnarray*}

H$\ddot{o}$lder's inequality after absorbing the dissipative term gives us the bound by

\begin{eqnarray*}
&&X(0) + c\int_{0}^{T}\lVert\partial_{3}u_{3}\rVert_{L^{s}}^{\frac{2}{\gamma - 2}}Xdt +c\int_{0}^{T}Z^{\frac{\gamma - 1}{\gamma - 2}}\lVert \nabla b\rVert_{L^{2}}^{2}dt\\
&&+c\int_{0}^{T}\lVert \partial_{3}u_{3}\rVert_{L^{s}}^{\frac{2}{\gamma - 2}}Xdt[\left(\int_{0}^{T}\lVert\Delta u\rVert_{L^{2}}^{2}dt\right)^{\frac{1}{4}} +\left(\int_{0}^{T}\lVert\Delta b\rVert_{L^{2}}^{2}dt\right)^{\frac{1}{4}}]\\
&&+c\left(\int_{0}^{T}Z^{\frac{\gamma - 1}{\gamma -2}}Xdt\right)^{\frac{\gamma - 2}{\gamma - 1}} \left(\int_{0}^{T}\lVert\Delta b\rVert_{L^{2}}^{2}dt\right)^{\frac{1}{\gamma - 1}}[\left(\int_{0}^{T}\lVert\Delta u\rVert_{L^{2}}^{2}dt\right)^{\frac{1}{4}} +\left(\int_{0}^{T}\lVert\Delta b\rVert_{L^{2}}^{2}dt\right)^{\frac{1}{4}}]
\end{eqnarray*}

We expand and bound the last term using Young's inequalities by

\begin{eqnarray*}
&&c\left(\int_{0}^{T}Z^{\frac{\gamma - 1}{\gamma - 2}}Xdt\right)^{\frac{\gamma - 2}{\gamma - 1}}\left(\int_{0}^{T}\lVert\Delta b\rVert_{L^{2}}^{2}dt\right)^{\frac{1}{\gamma - 1}}\left(\int_{0}^{T}\lVert\Delta u\rVert_{L^{2}}^{2}dt\right)^{\frac{1}{4}}\\
&&+c\left(\int_{0}^{T}Z^{\frac{\gamma - 1}{\gamma - 2}}Xdt\right)^{\frac{\gamma - 2}{\gamma - 1}}\left(\int_{0}^{T}\lVert\Delta b\rVert_{L^{2}}^{2}dt\right)^{\frac{1}{\gamma - 1}}\left(\int_{0}^{T}\lVert\Delta b\rVert_{L^{2}}^{2}dt\right)^{\frac{1}{4}}\\
&\leq&
c\left(\int_{0}^{T}Z^{\frac{\gamma - 1}{\gamma - 2}}Xdt\right)^{\frac{4(\gamma - 2)}{3\gamma - 7}} +\epsilon\int_{0}^{T}(\lVert \Delta u\rVert_{L^{2}}^{2} + \lVert\Delta b\rVert_{L^{2}}^{2})dt
\end{eqnarray*} 

Thus, after absorbing the dissipative term and applying Young's inequality on the middle term we obtain

\begin{eqnarray*}
&&X(T) + \int_{0}^{T}\lVert \Delta u\rVert_{L^{2}}^{2} + \lVert \Delta b\rVert_{L^{2}}^{2}dt\\
&\leq& X(0) + c\int_{0}^{T}\lVert\partial_{3}u_{3}\rVert_{L^{s}}^{\frac{2}{\gamma - 2}}Xdt + c\int_{0}^{T}Z^{\frac{\gamma - 1}{\gamma - 2}}\lVert \nabla b\rVert_{L^{2}}^{2}dt\\
&&+c\left(\int_{0}^{T}\lVert \partial_{3}u_{3}\rVert_{L^{s}}^{\frac{2}{\gamma - 2}}Xdt\right)^{\frac{4}{3}} + \epsilon(\int_{0}^{T}\lVert\Delta u\rVert_{L^{2}}^{2} + \lVert\Delta b\rVert_{L^{2}}^{2}dt) + c\left(\int_{0}^{T}Z^{\frac{\gamma - 1}{\gamma - 2}}Xdt\right)^{\frac{4(\gamma - 2)}{3\gamma - 7}}
\end{eqnarray*}

On the fourth and sixth terms, we use H$\ddot{o}$lder's inequalities and (8) to finally obtain

\begin{eqnarray*}
X(T) + \int_{0}^{T}\lVert \Delta u\rVert_{L^{2}}^{2} + \lVert \Delta b\rVert_{L^{2}}^{2}dt &\leq& X(0) + c\int_{0}^{T}\lVert\partial_{3}u_{3}\rVert_{L^{s}}^{\frac{2}{\gamma - 2}}Xdt +c\int_{0}^{T}Z^{\frac{\gamma - 1}{\gamma - 2}}\lVert \nabla b\rVert_{L^{2}}^{2}dt\\
&&+c\int_{0}^{T}\lVert \partial_{3}u_{3}\rVert_{L^{s}}^{\frac{8}{3}(\frac{1}{\gamma - 2})}Xdt +c\int_{0}^{T}Z^{\frac{4(\gamma - 1)}{3\gamma - 7}}Xdt
\end{eqnarray*}

Lemma 2.3 implies that if we take $\gamma > \frac{7}{3}$,  

\begin{equation*}
\sum_{i=2}^{3}\int_{0}^{T}\lVert \partial_{i}u_{i}\rVert_{L^{s}}^{\frac{8}{3\gamma - 7}}dt < \infty
\end{equation*}

suffices to complete the proof with Gronwall's inequality. 

\subsection{Proof of Corollary 1.3}
Corollary 1.3 is immediate from the following special case of the lemma due to [16]:

\begin{lemma}
Let u be a divergence-free sufficiently smooth vector field in $\mathbb{R}^{3}$. Then there exists a constant $C = C(q)$ such that 

\begin{equation*}
\lVert \partial_{i}u_{j}\rVert_{L^{q}} \leq C(\lVert w_{3}\rVert_{L^{q}} + \lVert \partial_{3}u_{3}\rVert_{L^{q}})
\end{equation*}

for $1 < q < \infty, i, j = 1, 2$. 

\end{lemma}

Thus, 

\begin{equation*}
\lVert \partial_{2}u_{2}\rVert_{L^{p}}^{r} \leq c(\lVert w_{3}\rVert_{L^{p}} + \lVert \partial_{3}u_{3}\rVert_{L^{p}})^{r} \leq c(\lVert w_{3}\rVert_{L^{p}}^{r} + \lVert \partial_{3}u_{3}\rVert_{L^{p}}^{r})
\end{equation*}

by Lemma 2.3 as $r \geq 0$. 

\subsection{Proof of Theorem 1.4}

Taking an $L^{2}$-inner product with $u$ on (2), integrating in time we see that 

\begin{equation}
\sup_{t\in [0,T]}\lVert u(t)\rVert_{L^{2}}^{2} + \int_{0}^{T}\lVert \Lambda^{\alpha}u\rVert_{L^{2}}^{2}dt < c(u_{0})
\end{equation}

\subsubsection{Estimate of $\lVert\nabla_{h}u\rVert_{L^{2}}^{2}$}

Taking an inner product with $-\Delta_{h}u$ on (2), applying Lemma 2.1 and integration by parts just like $J_{1}$ estimate in the previous proofs, Lemma 2.2 give us

\begin{eqnarray*}
&&\frac{1}{2}\partial_{t}\lVert \nabla_{h}u\rVert_{L^{2}}^{2} + \nu\lVert \Lambda^{\alpha}\nabla_{h}u\rVert_{L^{2}}^{2}\\
&\leq& c\lVert\partial_{3}u_{3}\rVert_{L^{s}}^{\frac{1}{\gamma}}\lVert\nabla u \rVert_{L^{2}}^{\frac{\gamma - 2}{\gamma}}\lVert \nabla_{h}u\rVert_{L^{2}}^{(1-\frac{1}{\alpha})(1+\frac{2}{\gamma})}\lVert\Lambda^{\alpha}\nabla_{h}u\rVert_{L^{2}}^{(\frac{1}{\alpha})(1+\frac{2}{\gamma})}
\end{eqnarray*}

due to Gagliardo-Nirenberg inequality. Young's and Gronwall's inequalities give us simiarly as before 

\begin{eqnarray*}
\sup_{t\in [0,T]}\lVert \nabla_{h}u(t)\rVert_{L^{2}}^{2} + \nu\int_{0}^{T}\lVert\Lambda^{\alpha}\nabla_{h}u\rVert_{L^{2}}^{2}dt \leq c(T) + c(T)\int_{0}^{T}\lVert \partial_{3}u_{3}\rVert_{L^{s}}^{\frac{2}{\gamma - 2}}\lVert \nabla u\rVert_{L^{2}}^{2}dt
\end{eqnarray*}

\subsubsection{Estimate on $\lVert\nabla u\rVert_{L^{2}}^{2}$}

Next, taking $L^{2}$-inner product on (2) with $-\Delta u$, 

\begin{equation*}
\frac{1}{2}\partial_{t}\lVert \nabla u\rVert_{L^{2}}^{2} + \nu\lVert \Lambda^{\alpha}\nabla u\rVert_{L^{2}}^{2} \leq c\int \lvert u_{3}\rvert \lvert \nabla u\rvert \lvert \nabla\nabla_{h}u\rvert + \lvert \nabla_{h}u\rvert \lvert \partial_{3}u\rvert^{2} = J_{1} + J_{2}
\end{equation*}

where 

\begin{eqnarray*}
J_{2} = c\int \lvert \nabla_{h}u\rvert \lvert \partial_{3}u\rvert^{2} \leq c\lVert \nabla_{h}u\rVert_{L^{2}}\lVert \nabla u\rVert_{L^{4}}^{2} \leq c\lVert \nabla_{h}u\rVert_{L^{2}}\lVert \nabla u\rVert_{L^{\frac{6}{5-2\alpha}}}^{(\frac{4\alpha - 3}{4})2}\lVert \Lambda^{\alpha}u\rVert_{L^{6}}^{(\frac{7-4\alpha}{4})2}
\end{eqnarray*}

by H$\ddot{o}$lder's and Gagliardo-Nirenberg inequality. Now Sobolev embedding and (7) combined with previous estimate give

\begin{eqnarray*}
\frac{1}{2}\partial_{t}\lVert \nabla u\rVert_{L^{2}}^{2} + \lVert \Lambda^{\alpha}\nabla u\rVert_{L^{2}}^{2} &\leq& c(\lVert \partial_{3}u_{3}\rVert_{L^{s}}^{\frac{2}{\gamma - 2}}\lVert \nabla u\rVert_{L^{2}}^{2} + \lVert \nabla_{h}u\rVert_{L^{2}}^{2})\\
&+& c\lVert \nabla_{h}u\rVert_{L^{2}}\lVert\Lambda^{\alpha}u\rVert_{L^{2}}^{(\frac{4\alpha - 3}{4})2}\lVert \nabla_{h}\Lambda^{\alpha}u\rVert_{L^{2}}^{\frac{4}{3}(\frac{7-4\alpha}{4})}\lVert\nabla\Lambda^{\alpha}u\rVert_{L^{2}}^{\frac{2}{3}(\frac{7-4\alpha}{4})}
\end{eqnarray*}

We integrate in time to obtain

\begin{eqnarray*}
&&\lVert \nabla u(T)\rVert_{L^{2}}^{2} + \int_{0}^{T}\lVert \Lambda^{\alpha}\nabla u\rVert_{L^{2}}^{2}dt\\
&\leq& c\int_{0}^{T}\lVert\partial_{3}u_{3}\rVert_{L^{s}}^{\frac{2}{\gamma -2}}\lVert \nabla u\rVert_{L^{2}}^{2}+ \lVert \nabla u\rVert_{L^{2}}^{2} dt \\
&&+c\sup_{t\in [0,T]}\lVert \nabla_{h}u\rVert_{L^{2}}\left(\int_{0}^{T}\lVert \Lambda^{\alpha}u\rVert_{L^{2}}^{(\frac{4\alpha - 3}{4})2(\frac{4}{4\alpha - 3})}dt\right)^{\frac{4\alpha - 3}{4}}\\
&&\times \left(\int_{0}^{T}\lVert\nabla_{h}\Lambda^{\alpha}u\rVert_{L^{2}}^{\frac{4}{3}(\frac{7-4\alpha}{4})(\frac{6}{7-4\alpha})}dt\right)^{\frac{7-4\alpha}{6}}\left(\int_{0}^{T}\lVert\nabla\Lambda^{\alpha}u\rVert_{L^{2}}^{\frac{2}{3}(\frac{7-4\alpha}{4})(\frac{12}{7-4\alpha})}dt\right)^{\frac{7-4\alpha}{12}}
\end{eqnarray*}

by H$\ddot{o}$lder's inequality. By the previous estimate and (9) we have the second term bounded by

\begin{eqnarray*}
&&c\left(\int_{0}^{T}\lVert \partial_{3}u_{3}\rVert_{L^{s}}^{\frac{2}{\gamma - 2}}\lVert \nabla u\rVert_{L^{2}}^{2}dt\right)^{\frac{5-2\alpha}{3}} \left(\int_{0}^{T}\lVert \nabla \Lambda^{\alpha}u\rVert_{L^{2}}^{2}dt\right)^{\frac{7-4\alpha}{12}}\\
&\leq& c\left(\int_{0}^{T}\lVert \partial_{3}u_{3}\rVert_{L^{s}}^{\frac{2}{\gamma - 2}}\lVert \nabla u\rVert_{L^{2}}^{2}dt\right)^{\frac{4(5-2\alpha)}{5+4\alpha}} + \left(\int_{0}^{T}\lVert\nabla \Lambda^{\alpha}u\rVert_{L^{2}}^{2}dt\right)
\end{eqnarray*}

by Young's inequality. Absorbing the last term, we have the bound of

\begin{eqnarray*}
&&c\int_{0}^{T}\lVert \partial_{3}u_{3}\rVert_{L^{s}}^{\frac{2}{\gamma - 2}}\lVert \nabla u\rVert_{L^{2}}^{2} + \lVert \nabla u\rVert_{L^{2}}^{2}dt\\
&&+c\left(\int_{0}^{T}\lVert \partial_{3}u_{3}\rVert_{L^{s}}^{\frac{8}{\gamma - 2}(\frac{5-2\alpha}{5+4\alpha})}\lVert \nabla u\rVert_{L^{2}}^{2}dt\right)\left(\int_{0}^{T}\lVert \nabla u\rVert_{L^{2}}^{2}dt\right)^{\frac{15-12\alpha}{5+4\alpha}}
\end{eqnarray*}

by H$\ddot{o}$lder's inequality. Therefore,

\begin{eqnarray*}
&&\lVert \nabla u(T)\rVert_{L^{2}}^{2} + \int_{0}^{T}\lVert \Lambda^{\alpha}\nabla u\rVert_{L^{2}}^{2}dt\\
&\leq& c\int_{0}^{T}\lVert \partial_{3}u_{3}\rVert_{L^{s}}^{\frac{2}{\gamma - 2}}\lVert \nabla u\rVert_{L^{2}}^{2} + \lVert \nabla u\rVert_{L^{2}}^{2}dt +c\left(\int_{0}^{T}\lVert \partial_{3}u_{3}\rVert_{L^{s}}^{\frac{8}{\gamma - 2}(\frac{5-2\alpha}{5+4\alpha})}\lVert \nabla u\rVert_{L^{2}}^{2}dt\right)
\end{eqnarray*}

by the Gagliardo-Nirenberg inequality and H$\ddot{o}$lder's inequality and (9). Gronwall's inequality and (5) complete the proof of Theorem 1.4.

\end{document}